\documentclass[]{article}
\usepackage{amsmath,amssymb,amsfonts,amsthm,bbm}
\numberwithin{equation}{section}
\usepackage[linesnumbered,ruled,algosection]{algorithm2e}
\usepackage[font=small,labelfont=bf]{caption}
\usepackage{array}
\usepackage{float}
\usepackage{graphicx}
\usepackage{epstopdf}
\usepackage{color,soul}
\usepackage{listings}
\usepackage[dvipsnames]{xcolor}

\definecolor{dkgreen}{rgb}{0,0.6,0}
\definecolor{gray}{rgb}{0.5,0.5,0.5}
\definecolor{mauve}{rgb}{0.58,0,0.82}
\def\EE{{\mathbb E}}
\def\PP{{\mathbb P}}

\def\RR{{\mathbb R}}
\def\rr{\mathbb R}
\def\nn{{\mathbb N}}
\def\cf{{\cal F}}

\begin{document}

\title{Optimal Forward Trading and Battery Control Under Renewable Electricity Generation}
\author{Juri Hinz \\ Juri.Hinz@uts.edu.au \and Jeremy Yee \\ jeremyyee@outlook.com.au}
\date{\today}

\maketitle

\begin{abstract}
    The increased market penetration of renewable energy sources and
    the rapid development of electric battery storage technologies
    yield a potential for reducing electricity price volatility while
    maintaining stability of the power grid.  This work presents an
    algorithmic approach to control battery levels and forward
    positions to optimally manage power output fluctuations caused by
    intermittent renewable energy generation. This paper will also
    explore the effect of battery technology on the firm's optimal
    trading behaviour in the electricity spot market.
 \end{abstract}

\section{Introduction} 

The recent proliferation of renewable energies and technological
progress in electric battery storage systems create an increasing
demand for sound algorithmic solutions to optimal control problems
arising in the dispatch optimization of power supply given a storage
facility and uncertainty caused by the market prices or/and weather
conditions.  Such problems are numerically challenging due to high
dimensionality of the state spaces involved. This work suggests a
quantitative approach to address a growing need for efficient
decentralized electricity dispatch and storage. 

\begin{figure}[h]
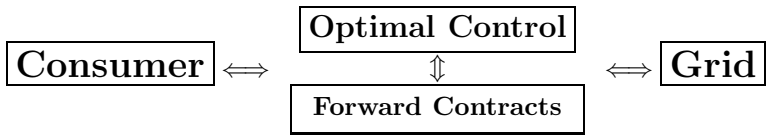

	$$
	\boxed{\text{{\bf \Large Consumer}}} 
	\Longleftrightarrow \begin{array}{c}
	\boxed{\text{\bf \large Optimal Control}}          \\
	\Updownarrow          \\
	\boxed{ \begin{array}{c} {\text{\bf Forward Contracts}} \end{array}}
	\end{array}
	\Longleftrightarrow
	\boxed{\hbox{\bf \Large Grid}}
	$$
	\caption{Traditional energy dispatch.}
        \label{old_grid}
\end{figure}

\begin{figure}[h]
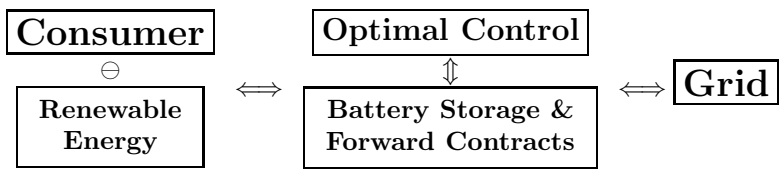

$$ 
{\begin{array}{c}
   \boxed{\text{{\bf \Large Consumer}}} \\ 
   \ominus\\
   \boxed{\begin{array}{c}
            \text{{\bf Renewable}} \\
            \text{{\bf Energy}}
          \end{array}
   }
 \end{array}}
\Longleftrightarrow \begin{array}{c}
\boxed{\text{\bf \large Optimal Control}} \\
\Updownarrow \\
\boxed{ \begin{array}{c} {\text{\bf Battery Storage \&}} \\ {\text{\bf Forward Contracts}} \end{array}}
\end{array}
\Longleftrightarrow   
\boxed{\hbox{\bf \Large Grid}}
$$
\caption{Energy dispatch in the presence of renewable energy and battery storage.}
\label{new_grid}
\end{figure}

Let us describe its typical framework. The traditional electricity
market players satisfy consumers' energy demand by purchasing
electricity in advance, usually taking positions in the so-called
day-ahead market (also called the spot market) such that any energy
imbalances must be compensated in real-time as they occur. This
real-time balancing can either be achieved through complex
over-the-counter trading or, more realistically, by transferring
supply from or to electricity grid at the so-called real time grid
prices. Figure \ref{old_grid} provides a simplified illustration of
this optimal control problem.  However, in the presence of storage and
renewable generation facilities, the problem changes. On this new
structure, the agent's control problem now requires simultaneously
taking optimal positions and setting optimum energy storage levels as
shown in Figure \ref{new_grid}. The decision optimization problem
becomes significantly more complex due to the uncertainty stemming
from the future battery capacity levels, electricity prices, and
output of renewable energy.

Many renewable energy sources such as wind and solar are notoriously
intermittent and unreliable. The potential of energy storage devices
to address the highly erratic nature of renewable energy generation
\cite{Breton_Moe_2009,Dincer_2011} and energy demand has been
discussed extensively in the literature (see
\cite{Beaudin_etal_2010,DiazGonzalez_etal_2012,Evans_etal_2012,Irena_2015,
  Yang_etal_2014}). Their incorporation into a modern energy grid will
encourage more environmentally friendly policies which will also have
significant impact on investor atittudes towards firms
\cite{Ramiah_etal_2013, Chan_Walter_2014,Renneboog_etal_2008}.  The
authors of \cite{Lu_etal_2014} studied the possible usage of battery
storage systems to defer costly modifications to the energy grid by
addressing peak loads in the power grid. An extensive recent review of
available energy storage technologies has been given by
\cite{Luo_etal_2015} and future innovation looks bright. While there
exist numerous types of energy storage systems,
\cite{Beaudin_etal_2010} found that no single storage system
consistently outperforms all the others for all types of renewable
energy sources and applications. So for the sake of simplicity, this
paper will assume that the energy retailer pictured in Figure
\ref{new_grid} uses a battery device for storing energy.  However, the
methods and results contained within this paper can easily be extended
for other types of storage technologies or even to the use of multiple
types of storage devices. From a real options analysis point of view,
the incorporation of energy storage devices into energy grid also
poses interesting investment questions. The work done by
\cite{Bakke_etal_2016,Bradbury_etal_2014,Locatelli_etal_2016} examined
the profitability of investing in energy storage devices.  However,
\cite{Schachter_Mancarella_2016} questions the suitability of the
current real options approach, stating that the risk neutrality
assumption may not be appropriate for risk averse investors. The
introduction of batteries also gives rise to important optimal
stochastic control problems. The optimal dynamic storage and discharge
of energy from battery devices has been examined in
\cite{Dokuchaev_2016,Kim_Powell_2011,Oudalov_etal_2007,Teleke_etal_2010}.

Rather than focusing on capacity investment decision, the present
contribution focuses on optimal operational management in terms of
energy purchase and dispatch optimization, given a storage device of a
fixed capacity. Thereby, we suppose that storage facility is only used
for the compensation of any imbalance between consumers' demand,
renewable energy generation and an existing financial position. This
issue is connected to the market behaviour addressed in
\cite{Cartea_Villaplana_2008} in terms of reducing the risk in the
sense of \cite{Benth_etal_2008}, since a storage acts as a safety
buffer. The present work helps investigating the effect of battery on
forward energy trading.

This paper is organized as follows: Section \ref{sec:problem}
introduces the model while Section \ref{sec:markov} frames the main
question as a Markov Decision Problem whose numerical solution is
outlined in Section \ref{sec:numerics}. Section
\ref{sec:illustration} provides a numerical study of prices which
exhibit distinct mean reversion and seasonality.  Finally, Section
\ref{sec:discussion} examines the impact of battery on the optimal
forward energy trading with conclusions provided in Section
\ref{sec:conclusion}.

\section{Problem Setting} \label{sec:problem}

Within a given time horizon $t=0, \dots, T-1$, the net energy demand
$Q_t$ within each period is the difference between the consumer's
demand and the renewable energy output. Given an existing financial
position $F_t$, the energy imbalance $F_t-Q_t$ will be compensated
using energy from the battery storage followed by a possible offset
through real-time energy from the power grid.  That is, in the case of
energy surplus $F_t-Q_t \ge 0$, the electricity is first used to
charge the battery up to the maximal level and the remaining energy is
then supplied to the grid. Similarly, if there is an energy shortage
$F_t-Q_t <0$, the required electricity is taken from the battery up to
a minimal battery level before the required electricity rest is taken
from the grid.

Let us assume that the net demand realization is given by
$Q_t=q_t+\varepsilon_t$ with zero-mean random variable $\varepsilon_t$
describing the deviation of the net demand from its predicted level
$q_t$ and suppose that the financial position is given in terms of
$F_t=q_t + l$ where the quantity $l$ describes a safety margin and
stands retailer's decision to buy/sell in the spot market an energy
amount $q_t + l$ which deviates from the predicted net demand $q_t$ by
$l$.  Thereby, we model the decision of the retailer in the choice
$F_t=F_t(a)$ of financial positions in terms of the action
$a \in \mathbf{A}$ from a finite set $\mathbf{A}$ of all possible
actions, each characterized by its specific safety margin $l(a)$.
With the assumptions above, given the action $a \in \mathbf{A}$, the
realized net energy to be balanced is given by
$$
F_t(a)-Q_t= q_t+l(a)-(q_t + \varepsilon_t )=l(a) - \varepsilon_t, 
\quad a \in \mathbf{A}.
$$
That is, the action $a \in \mathbf{A}$ determines a certain
distribution $\nu_t(a)$ of the energy volume which must be balanced
and is determined as $\nu_t(a) \sim l(a)-\varepsilon_t$ for all
$a \in \mathbf{A}$.  In order to describe the battery storage control,
we suggest discretizing the storage levels by a finite set
$\mathbf{P}$. Having chosen the action $a \in \mathbf{A}$, the
imbalance energy $F_t(a)-Q_t=l(a) - \varepsilon_t$ follows a
distribution $\nu_t(a)$ which determines for each battery storage
level $p \in \mathbf{P}$ a probability $\alpha_{p, p'}^a$ that the
storage reaches its next-day level $p' \in \mathbf{P}$. Furthermore,
the expected energy excess $\underline e^a_{p}$ and shortage
$\overline e^a_{p}$ are uniquely determined by the current battery
level $p \in \mathbf{P}$ and the action $a \in \mathbf{A}$ through the
imbalance distribution $\nu_t(a)$.
 
Now let us turn to the costs of energy imbalance. For this, we
introduce the random variables
$$0 \le \underline \Pi_t \le \overline \Pi_t, \qquad t=0, \dots, T$$
which stand for the sell/buy real time grid prices expected at time
$t$, when the financial position $F_t(a)$ is taken.  With these
definitions, the revenue/costs associated with energy imbalance for
the action $a \in \mathbf{A}$ are modeled by
$ -\overline e^a_{p} \overline \Pi_t + \underline e^a_{p} \underline
\Pi_t.  $ Finally, let us denote by  $\Pi_t$ the energy price at time $t=0, \dots, T$.
 Since we
assume that all feasible financial positions are given as $q_t+l(a)$
with $a \in \mathbf{A}$, the position costs for the action
$a \in \mathbf{A}$ are
$$
(q_t+l(a)) \Pi_t=q_t \Pi_t+l(a)\Pi_t.
$$
With assumptions and notations as above, the revenue/loss associated
with the action $a \in \mathbf{A}$ depends on the current  price
$\Pi_t$, the expected demand $q_t$, and the recent battery level
$p \in \mathbf{P}$ as
\begin{equation}
r_t(p, (q_t, \Pi_t), a)= -q_t  \Pi_t   - l(a) \Pi_t - \overline e^a_{p} \overline \Pi_{t}  + 
\underline e^a_{p} \underline \Pi_{t}. \nonumber
\end{equation}
Observe that the term $- q_t \Pi_t $ neither depends on the action $a$
nor on battery level $p$. Thus, we agree that the choice
$a \in \mathbf{A}$ of the optimal safety margin $l(a)$ will depend
only on electricity price and re-define the reward as
\begin{equation}
r_t(p, \Pi_t, a)=- l(a) \Pi_t - \overline e^a_{p} \overline \Pi_{t}  + 
\underline e^a_{p} \underline \Pi_{t}. \label{rewfct} 
\end{equation}
Note that we also do not consider other revenues associated with income streams 
due to delivery commitments at fixed price. On this account, the revenue (\ref{rewfct})
serves a vehicle to optimize trading activity in terms of optimal safety margins 
and does not reflect actual  cash flows.  

The revenue optimization from battery storage management is a typical
sequential decision problem under uncertainty.  Having chosen at time
$t=0, \dots, T-1$ an action $a \in \mathbf{A}$ in the data 
$(p,  \Pi_t)$ a certain revenue/costs $r_t(p,  \Pi_t, a)$
is incurred immediately. However, the action $a \in \mathbf{A}$ also
changes the probability of transition to the subsequent states (next
battery levels) which influences all future revenues and
decisions. Problems of this type are naturally formulated and solved
in terms of the so-called Markov Decision Theory.  In what follows, we
formulate our storage control problem within this standard framework.

\section{Markov Decision Theory for Battery Control}
\label{sec:markov}

Let us review the classical finite-horizon Markov decision theory
following \cite{BRieder}.  On a finite time horizon $0,\dots,T$,
consider a random dynamics whose state $x$ evolves in (measure) space
$E$ and is controlled by actions $a$ from a finite action set
$\mathbf{A}$.  For each $a \in \mathbf{A}$, we assume that
$K^{a}_{t}(x, dx')$ is a stochastic transition kernel on $E$. A
mapping $\pi_{t}:E \mapsto \mathbf{A}$ which describes the action that
the controller takes at time $t$ is called a \emph{decision rule}. A
sequence of decision rules $\pi=(\pi_{t})_{t =0}^{T-1}$ is called a
\emph{policy}. For each initial point $x_{0} \in E$ and each policy
$\pi$, there exists a probability measure
$\PP^{x_{0}, \pi}$ and a stochastic process $(X_{t})_{t = 0}^{T}$ such
that $\PP^{x_{0},\pi}(X_{0}=x_{0})=1$ and
\begin{equation}
\label{fo:chain}
\qquad \PP^{x_{0}, \pi}(X_{t+1} \in B \, | \, X_{0}, \dots, X_{t})=K^{\pi_{t}(X_{t})}_{t}(X_{t},B )
\end{equation}
holds for each (measurable) $B \subset E$ at all times
$t=0,\dots, T-1$. That is, given the system state $X_t$ at time $t$,
the action $a = \pi_{t}(X_{t})$ is used to pick the transition
probability $K_{t}^{a = \pi_{t}(X_{t})}(X_{t}, \cdot)$ which randomly
drives the system from $X_{t}$ to $X_{t+1}$ with the distribution
$K^{\pi_{t}(X_{t})}_{t}(X_t,\,\cdot\,)$.  Let us use
${\cal K}^{a}_{t}$ to denote the one-step transition operator
associated with the transition kernel $K^{a}_{t}$ when the action
$a\in \mathbf{A}$ is chosen.  In other words, for each action
$a \in \mathbf{A}$ the operator ${\cal K}^{a}_{t}$ acts on functions
$v$ by
\begin{equation}
({\cal K}^{a}_{t} v )(x)=\int_{E} v(x')K^{a}_{t}(x,dx') \qquad x \in E,\label{transop}
\end{equation}
whenever the above integrals are well-defined.  Now, let us turn to
the definition of the control costs.  For each time $t$, we are given
the \emph{$t$-step reward function}
$r_{t}: E \times \mathbf{A} \mapsto \RR$, where $r_{t}(x, a)$
represents the reward for applying an action $a \in \mathbf{A}$ when
the state of the system is $x \in E$ at time $t$.  At the end of the
time horizon, at time $T$, it is assumed that no action can be
taken. Here, if the system is in a state $x$, a {\it scrap value}
$r_{T}(x)$, which is described by a pre-specified {\it scrap function}
$r_{T}: E \to \RR$, is collected.

Given an initial point $x_{0}$, the goal is to maximize the expected
finite-horizon total reward
\begin{equation}
  \label{policyvalue}
  v_0^\pi(x_0) =\EE^{x_{0}, \pi}\left(\sum_{t=0}^{T-1} r_{t}(X_{t}, \pi_{t}(X_{t})) + r_{T}(X_{T}) \right),
\end{equation}
over all possible policies $\pi=(\pi_t)_{t=0}^{T-1}$, where
$\EE^{x_{0},\pi}$ denotes the expectation over the controlled Markov
chain defined by (\ref{fo:chain}).  In other words, to find the
argument $\pi^{*}=(\pi^{*}_{t})_{t=0}^{T-1}$ such that
\begin{equation}
  \label{rewmax}
  \pi^*(x_0)=\arg\max_{\pi}v_0^\pi(x_0) , \quad x_0 \in E.
\end{equation}
The maximization (\ref{rewmax}) is well-defined under additional
assumptions (see \cite{BRieder}, p. 199).

The calculation of the optimal policy is addressed in the following
setting.  For $t=0, \dots, T-1$, introduce the {\it Bellman operator}
\begin{eqnarray}
{\cal T}_{t}v(x)=\sup_{a \in \mathbf{A}} \left(r_{t}(x,a)+ {\cal K}^{a}_{t}v(x) \right), \quad x \in E
\end{eqnarray}
which acts on each measurable function $v: E \to \rr$ where the
integrals ${\cal K}_{t}^{a}v$ for all $a \in \mathbf{A}$ exist.
Further, consider the {\it Bellman recursion}
\begin{eqnarray}
v^{*}_{T}=r_{T},  \label{bel1} \quad
v^{*}_{t}= {\cal T}_{t} v^{*}_{t+1} \qquad \hbox{for $t=T-1, \dots, 0$.}
\end{eqnarray}
Under appropriate assumptions, there exists a recursive solution
$(v^{*}_{t})_{t=0}^{T}$ to the Bellman recursion, which gives the
so-called {\it value functions} and determines an optimal policy
$\pi^{*}$ via
$$
\pi^{*}_{t}(x)={\rm argmax}_{a \in \mathbf{A}}\left(r_{t}(x,a)+ {\cal K}_{t}^{a}v^{*}_{t+1}(x) \right), \quad x \in E
$$
for all $t=0, \dots, T-1$.

Consider now a Markov decision model whose state evolution consists of
one discrete and one continuous component. To be more specific, we
assume that the state space $E=\mathbf{P} \times \rr^{d}$ is the
product of a finite space $\mathbf{P}$ and the Euclidean space
$\rr^{d}$.  We suppose that the discrete component $p \in \mathbf{P}$
is driven by a finite number of actions $a \in \mathbf{A}$ in terms
stochastic matrices
$$
(\alpha_{p, p'}^{a})_{p, p' \in \mathbf{P}}, \quad a \in \mathbf{A}
$$
where $\alpha_{p, p'}\in [0, 1]$ stands for the transition probability
from $p \in \mathbf{P}$ to $p' \in \mathbf{P}$ if the action
$a\in \mathbf{A}$ was taken.  Furthermore, we assume that the
continuous state component evolves as an uncontrolled Markov process
$(Z_{t})_{t=0}^{T}$ on $\rr^{d}$ realized on a probability space
$(\Omega, \cf, \PP)$ whose evolution is driven by random linear
transformations
\begin{align}
Z_{t+1}=W^{}_{t+1}Z_{t} \label{lintrafo}
\end{align}
with pre-specified independent and integrable disturbance matrices
$(W^{}_{t})_{t=1}^{T}$.  In this setting, the transition and the
Bellman operators are given by
\begin{align}
{\cal K}_{t}^{a}v(p,z)&=\sum_{p' \in \mathbf{P}}\alpha_{p, p'}^{a}\EE_{\PP}(v(p', W^{}_{t+1}z)),  \quad z \in \rr^{d}  \nonumber \label{transop1} \\
{\cal T}_{t}v(p, z)&=\max_{a \in  \mathbf{A}}(r_{t}(p, z, a)+ \sum_{p' \in \mathbf{P}} \alpha_{p, p'}^{a}\EE_{\PP}(v(p', W^{}_{t+1}z)) ). \nonumber
\end{align}
for $ t=0, \dots, T-1$, and $a \in\mathbf{A}$.
 Finally, let us assume that the reward and scrap
functions
\begin{equation}
z \mapsto r_{t}(p, \cdot ,a), \enspace z \mapsto r_{T}(p, z) \quad t=0,\dots,T-1, \quad p\in \mathbf{P}, \quad a \in \mathbf{A}
\end{equation}
are convex and globally Lipschitz continuous in the continuous
component $z \in \rr^{d}$ of the state variable $(p, z)$.  Such Markov
decision problems are referred to as {\it convex switching systems}
(see \cite{HinzConvexSwitching}).  For such system, the backward
induction for $ p \in \mathbf{P}$, $z \in \rr^d$
\begin{align}
	v^{*}_{T}(p, z)&=r_{T}(p,z),  \label{use1} \\
	v^{*}_{t}(p,z)&=\max_{a \in \mathbf{A}}(r_{t}(p, z, a)+ \sum_{p' \in \mathbf{P}} \alpha_{p, p'}^{a}\EE_\PP(v^*_{t+1}(p', W^{}_{t+1}z)) ) \label{use2}
\end{align}
for $ t=T-1, \dots 0$ yields value functions $(v^*_t)_{t=0}^T$ which
provide an optimal policy $(\pi^*_t)_{t=0}^T$ via
\begin{equation}
\pi^*_t(p, z)={\rm argmax}_{a \in \mathbf{A}}(r_{t}(p, z, a)+ \sum_{p' \in \mathbf{P}} \alpha_{p, p'}^{a}\EE_\PP(v^*_{t+1}(p', W^{}_{t+1}z)) ).  \label{use3}
\end{equation}

\section{Numerical Solution and Diagnostics} \label{sec:numerics}

This paper will use the numerical approaches studied in
\cite{HinzConvexSwitching,HinzYapAlgorithms,HinzYeeUsersGuide} to
solve Markov Decision problems of \emph{convex switching} type
introduced in the previous section. We refer interested readers to
those works for a more detailed explanation. However, for the sake of
this paper's completeness, this section will briefly outline these
methods. The first step in obtaining a numerical solution to the
backward induction
\begin{eqnarray*}
	v^{*}_{T}=r_{T}, \qquad
	v^{*}_{t}={\cal T}_{t}v^{*}_{t+1}, \qquad t=T-1, \dots 0
\end{eqnarray*}
is an appropriate discretization of the Bellman operator
\begin{equation}
{\cal T}_{t}v(p, z)=\max_{a \in \mathbf{A}}\left(r_{t}(p, z, a)+ \sum_{p' \in \mathbf{P}} \alpha_{p, p'}^{a}\EE(v(p', W^{}_{t+1}z)) \right). \nonumber
\end{equation}
For this reason, we consider a modified Bellman operator
${\cal T}^{n}_{t}$ instead of ${\cal T}_{t}$ with the expectation
$\EE(v(p', W^{}_{t+1}z))$ replaced by its numerical counterpart as
$$
{\sum_{n=1}^{N} \nu^{N}_{t+1}(n) v(p', W^{}_{t+1}(n)z)} )
$$
defined in terms of an appropriate distribution sampling
\begin{equation}
\text{$(W^{}_{t+1}(n))_{n=1}^{N}$ of each disturbance $W_{t+1}$ with weights $(\nu^{N}_{t+1}(n))_{n=1}^{N}$.}  \label{distsamp}
\end{equation}
In the resulting modified  backward induction
\begin{eqnarray}
v^{}_{T}=r_{T}, \qquad
v^{}_{t}={\cal T}^{N}_{t}v^{}_{t+1}, \qquad t=T-1, \dots 0 \label{modval}
\end{eqnarray}
the functions $(v_{t})_{t=0}^{T}$ need to be described by
algorithmically tractable objects. We may then approximate these
convex functions in terms of piecewise linear and convex functions in
the following manner. First, we introduce the so-called sub-gradient
envelope ${\cal S}_{G}f$ of a convex function $f: \rr^d \to \rr$ on a
grid $G \subset \rr^d$ as
$$
{\cal S}_{G}f=\vee_{g \in G} (\triangledown_{g}f)
$$ 
which is a maximum of the sub-gradients $\triangledown_{g}f$ of $f$ on
all grid points $g \in G$.  Using sub-gradient envelope operator, we
define the {double-modified Bellman operator} as
$$
{\cal T}^{G, N}_{t}v(p, \cdot)
=
{{\cal S}_{G}}
{\cal T}^{N}_{t}v(p, \cdot),
$$
where the operator ${\cal S}_{G}$ stands for the sub-gradient
envelope on the grid $G$.  The corresponding backward induction
\begin{eqnarray}
v^{}_{T}(p, \cdot)&=&{\cal S}_{G}r_{T}(p, \cdot),  \quad  p \in \mathbf{P}\label{scheme1}\\
v^{}_{t}(p, \cdot)&=&{\cal T}^{G,N}_{t}v_{t+1}(p, \cdot), \quad  p \in \mathbf{P}, \, t=T-1, \dots 0. \label{scheme2}
\end{eqnarray}
yields the so-called double-modified value functions
$(v_{t})_{t=0}^{T}$ which enjoy excellent algorithmic
properties. Namely, since the functions $(v^{}_{t})_{t=0}^{T}$ are
piece-wise linear and convex, they can be expressed using matrix
representations. Note that any piecewise convex function $f$ can be
described by a matrix where each of the linear functionals is
represented by one of the matrix's rows.  To denote this relation, let
us agree on the following notation: Given a function $f$ and a matrix
$F$, we write $f \sim F$ whenever $f(z)=\max (Fz)$ holds for all
$z\in\rr^d$.  It turns out that the sub-gradient envelope operation
${\cal S}_{G}$ on a grid $G$ corresponds to a specific
row-rearrangement operator in the following sense
$$
f \sim F \quad \Rightarrow \quad {\cal S}_{G}f \sim \Upsilon_{G}[F]
$$
where the row-rearrangement $\Upsilon_{G}$ associated with
 grid $G=\{g^1, \dots, g^m\}\subset \rr^{d}$ acts on matrix $F$ with $d$ columns as
follows:
\begin{equation} \label{row-rearrangment}
(
{\Upsilon_{G} }F)_{i, \cdot}= L_{ {\rm argmax}(F g^{i}), \cdot} \qquad \hbox{for all $g^i \in\{g^1, \dots, g^m\}= G$}
\end{equation}

For piecewise convex functions, the result of maximization, summation,
and composition with linear mapping, {followed by sub-gradient
  envelope} can be obtained using their matrix representatives. More
precisely, if
$$
f_{1} \sim F_{1}, \quad f_{2} \sim F_{2}
$$
holds, it follows that
\begin{eqnarray*}
	{\cal S}_{G}(f_{1}+f_{2}) & \sim & \Upsilon_{G}(F_{1})+ \Upsilon_{G}(F_{2}) \\
	{\cal S}_{G}(f_{1} \vee f_{2}) & \sim & \Upsilon_{G}(F_{1} \sqcup F_{2}) \\
	{\cal S}_{G}(f_{1}( W_{t+1}(k) \cdot ) & \sim & \Upsilon_{G}(F_{1} W_{t+1}(k))
\end{eqnarray*}
where the operator ${\sqcup}$ stands for binding matrices by rows,
which yields a matrix whose rows contain all rows from each
participating matrix. Using these relations, it turns out that the
double-modified backward induction can be rewritten in terms of the
row-rearrangement operator $\Upsilon=\Upsilon_{G}$, binding operator
$\sqcup$ and summations, applied to matrix representatives of the
double-modified value functions. Let us describe the resulting
algorithm.

Given a finite grid $G \subset \rr^d$, implement the row-rearrangement
operator $\Upsilon=\Upsilon_{G}$ and the matrix binding operator
$\sqcup$.  Determine a distribution sampling $(W^{}_{t}(n))_{n=1}^{N}$
of each disturbance $W_{t}$ with the corresponding weights
$( \nu^{N}_{t}(n))_{n=1}^{N}$ for $t=1, \dots, T$. Given reward
functions $(r_{t})_{t=0}^{T-1}$ and scrap value $r_{T}$, determine the
matrix representative of their sub-gradient envelopes
$$
{\cal S}_{G}r_{t}(p, \cdot, a) \sim  R_{t}(p, a), \qquad {\cal S}_{G}r_{T}(p, \cdot) \sim R_{T}(p)
$$
for $t=0, \dots, T-1$, $p \in \mathbf{P}$ and $a \in
\mathbf{A}$. Introduce matrix representatives of each value function
$$
v_{t}(p, \cdot) \sim V_{t}(p) \quad \text{for $t=0, \dots, T$, $p \in \mathbf{P}$,}
$$
which are obtained via Algorithm \ref{algo1}.

\vspace{3mm}
\begin{minipage}{0.95\textwidth}
	\setlength{\algomargin}{-0.1em}
	\begin{algorithm}[H]  
		\caption{Value Function Approximation \label{algo1}}
		\For {
			$p \in {P}$
		}
		{
			$V_T(p)  \leftarrow R_T(p) \leftarrow {\cal S}_G r_{T}(p,.)$ 
			
			\For {$ a \in  A, t=0, \dots, T$}{	
				$R_t(p, a) \leftarrow  {\cal S}_G r_{t}(p,., a)$
			}
		}
		\For{$t \in \{T-1,\dots,0\}$}{
			\For {$p \in {P}$}{
				$ V^{E}_{t+1}(p) \leftarrow \sum_{n=1}^N \nu^N_{t+1}(n) \Upsilon \left[
				V_{t+1}(p)  W_{t+1}(n)  \right] $
			}
			\For {$p \in  {P}$}{
				$ V_{t}(p) \leftarrow \Upsilon \sqcup_{a\in \mathbf{A}} \left( 
				R_t(p,a) +
				\sum_{p'\in {P}}\alpha_{p,p'}^a V^E_{t+1}(p') \right) $ 
			}
		}
	\end{algorithm}
\end{minipage}
\vspace{3mm}

Having calculated matrix representatives $(V_{t})_{t=0}^{T}$, the
approximations $( v^{}_{t})_{t =0}^{T}$, $( v^{E}_{t})_{t =0}^{T}$ of
the value functions and their expectations are given by
\begin{eqnarray}
v_{t}(p, z)&=& \max(V_{t}(p)z), \label{appval} \\
v^{E}_{t}(p, z)&=&\max( \sum_{n=1}^{N} \nu_{t}(n)
{\Upsilon}[V_{t}(p) \cdot  W_{t}(n)]   z) \label{appval1}
\end{eqnarray}
for all $z \in \rr^{d}$, $t=1, \dots, T$, and $p \in \mathbf{P}$.
Furthermore, an approximately optimal strategy
$(\tilde \pi_{t})_{t=0}^{T-1}$ is obtained for $t=0, \dots, T-1$ as
\begin{equation}
\pi_{t}(p,z)=\arg\max_{a \in \mathbf{A}}\left(r_{t}(p,z, a) + \sum_{p' \in \mathbf{P}}\alpha^{a}_{p,p'} v^{E}_{t+1}(p',z)\right).
\label{apstrat}
\end{equation}

We utilize an adaptation of the duality techniques developed by
C. Rogers \cite{Rogers2007} (see also \cite{RogersMathFin},
\cite{HaughKogan}, \cite{Chen2007}), to assess the quality of our
numerical solution. In its original formulation, the duality approach
provides an upper bound estimate for the unknown value function. This
technique has further developed in the context of discrete-time,
giving a promising view on duality of stochastic control which was
achieved in terms of the so-called \emph{information relaxation dual}
that was pioneered in the seminal paper by Brown, Smith, and Sun
\cite{BrownSunSmith}.

Here, we follow to the diagnostics method described in
\cite{HinzYeeUsersGuide} whose proofs are found in
\cite{HinzYapAlgorithms}.  Suppose that a candidate
$(\pi_{t})_{t=0}^{T-1}$ for approximatively optimal policy is
given. To estimate its distance-to-optimality, we address the
performance gap
$[v^{ \pi}_{0}(p_{0}, z_{0}), v^{\pi^{*}}_{0}(p_{0},z_{0})]$ in policy
values (\ref{policyvalue}) at a given starting point $z_{0}=Z_{0}$.
For this, we construct random variables
$ \underline{v}^{\pi, \varphi}_{0}(p_{0},z_{0})$,
$\overline{v}^{ \varphi}_{0}(p_{0},z_{0})$ satisfying
\begin{eqnarray}
\EE(\underline v^{\pi, \varphi}_{0}(p_{0},z_{0}))= v^{ \pi}_{0}(p_{0}, z_{0}) \le
v^{\pi^{*}}_{0}(p_{0},z_{0}) \le \EE(\bar v^{\pi, \varphi}_{0}(p_{0},z_{0})). \label{empbound}
\end{eqnarray}
The calculation of the expectations
\begin{equation} \label{exp}
\EE({\underline v^{\pi, \varphi}_{0}(p_{0}, z_{0})}) \text{  and  } \EE({\bar v^{\pi, \varphi}_{0}(p_{0}, z_{0})})
\end{equation} 
is realized through an efficient recursive Monte-Carlo scheme, which
yields approximations to (\ref{exp}) along with appropriate confidence
intervals.

For a practical application of this bound estimation, we assume that
an approximate solution yields a candidate $(\pi_{t})_{t=0}^{T-1}$ for
an optimal strategy, as in (\ref{apstrat}) based on approximations of
the value and of the expected value functions as in (\ref{appval}) and
(\ref{appval1}).  Further, choose a path number $K$ and a nesting
number $I \in \nn$ to obtain for each $k=1, \dots, K$ and
$i=0, \dots, I$ independent realizations $(w^{i,k}_{t})_{t=0}^{T}$ of
the random variables $(W_{t})_{t=0}^{T=1}$ and define for
$k=1, \dots, K$ the state trajectories $(z_{t}^{k})_{t=0}^{T}$
recursively
\[z_{0}^{k}:=z_{0}, \quad z^{k}_{t+1}=w^{0,k}_{t+1}z^{k}_{t}, \qquad t=0, \dots, T-1.\]

Estimators for the bounds in (\ref{empbound}) can be obtained using
Algorithm \ref{algo2} below.

\vspace{3mm}
\begin{minipage}{0.95\textwidth}
	\setlength{\algomargin}{-0.1em}
	\begin{algorithm}[H]
		\caption{Solution Diagnostics \label{algo2}}
		
		\For {$k = 1,\dots, K$}{
			$z_0^k \leftarrow z_0$ \\
			\For {$t=0, \dots, T-$}{
				$z_{t+1}^k \leftarrow w_{t+1}^{0,k} z_t^k$\\
				$	\varphi^{k}_{t+1}(p,a) \leftarrow \sum_{p' \in \mathbf{P}} \alpha^{a}_{p, p'}(
				\frac{1}{I}\sum_{i=1}^{I}v_{t+1}(p',  w^{i,k}_{t+1}z^{k}_{t})  
				- 
				v_{t+1}(p',  z^{k}_{t+1}))$ \\
			}
		}
		\For {$k =1,\dots, K$}{
			\For {$p \in  {P}$}{
				
				$\bar v^{k}_{T}(p)  \leftarrow
				\underline v^{k}_{T}(p) \leftarrow r_{T}(p, z^{k}_{T})
				$
			}
			\For{$t \in \{T-1,\dots,0\}$}{
				\For {$p \in  {P}$}{
					$
					\bar v^{k}_{t}(p) \leftarrow  \max_{a \in \mathbf{A}}\big[r_{t}(p, z^{k}_{t},a)+\varphi^{k}_{t+1}(p, a)         
					+ \sum_{p' \in \mathbf{P}} \alpha^{a}_{p, p'}
					\bar v^{k}_{t+1}(p')\big]$ 
					\\
					$a_t^k \leftarrow  \pi_{t}(p, z^{k}_{t})$\\
					$
					\underline{v}^{k}_{t}(p) \leftarrow r_{t}(p, z^{k}_{t},  a_t^k)+
					\varphi^{k}_{t+1}(p, a^k_t)	+
					\sum_{p' \in \mathbf{P}} \alpha^{a_t^k}_{p, p'}
					\underline{v}^{k}_{t+1}(p')
					$
				}
			}
		}    		  
		Determine estimators $\frac{1}{K}\sum_{k=1}^{K} \bar {v}^{k}_{0}(p)$ and
		$\frac{1}{K}\sum_{k=1}^{K} \underline{v}^{k}_{0}(p)$	to  (\ref{exp}).
	\end{algorithm}
\end{minipage}
\vspace{3mm}

\noindent
{\bf Remark:}
It is important to note that the convexity and Lipschitz continuity of
the reward and scrap functions presented in (\ref{rewfct}) is
essential for strong convergence properties of Algorithm
\ref{algo1}. However, these assumptions are not required for the
results from Algorithm \ref{algo2} to be valid. Using an appropriate
embedding of the state vector $Z_t$, the state dynamics presented in
(\ref{lintrafo}) is flexible enough to encompass a wide range of state
evolution process including geometric Brownian motion, auto regression
of order one and GARCH-like features. In addition, the dynamics
presented in (\ref{lintrafo}) can be extended to include a more
general function specification.

\section{Battery Control for Auto-Regressive State Dynamics} \label{sec:illustration} 

As a demonstration, let us consider a model based on the
auto-regressive state dynamics.  To cover this process under our
framework, we introduce $Z_t = (Z^{(1)}_t, Z^{(2)}_t) = (1, Z^{(2)}_t)$
where the first component equals to one for $t=0, \dots, T$ and define
the linear state dynamics
\begin{align} \label{rec}
Z_{t+1}=\begin{bmatrix} 1 \\ Z^{(2)}_{t+1} \end{bmatrix}
= \begin{bmatrix} 1 & 0 \\ 
\mu + \sigma N_{t+1} & \phi \end{bmatrix}
\begin{bmatrix} 1 \\ Z^{(2)}_t \end{bmatrix}=W_{t+1}Z_t, \quad t=0, \dots, T-1.
\end{align}
with  constants $\mu \in \rr$, $\sigma \in \rr_+$ and
$\phi \in [0, 1] $, driven by independent standard normally
distributed random variables $(N_t)_{t=1}^T$.  Further, we assume that
the electricity price $(\Pi_t)_{t=0}^T$ is governed by the
function $f:\mathbb{N}_+ \times \mathbb{R} \mapsto \mathbb{R}$ applied
to the state process as
\begin{equation}
\Pi_t=f(t, Z^{(2)}_t), \qquad t=0, \dots, T. \label{expauto}
\end{equation}
In this work, we restrict ourselves to consider deterministic affine
linear functions $(f(t, \cdot))_{t=0}^T$ to appropriately describe any
seasonal pattern of the electricity  price, frequently observed in practice. 

To model the the consumer's energy demand $Q_t = q_t + \varepsilon_t$
realized at time $t$, we suppose that conditioned on the information
at time $t$, the deviation $\varepsilon_t$ of the realized demand from
its predicted value $q_t$ follows a centered normal distribution with
a given variance $\varsigma^2 \in \rr_+$.  To describe the evolution
of the battery storage levels, let us assume that a finite set
$\mathbf{P}$ describes the storage levels of the battery which are
equidistantly spaced between the minimal
$\underline{p}=\min \mathbf{P}$ and the maximal
$ \overline{p}= \max \mathbf{P}$ level with a step size $\Delta>0$ .  Furthermore, consider a
finite set $\mathbf{A}$ of actions along with the function
$l: \mathbf{A} \to \rr$ prescribing the safety margin $l(a)$ chosen by
the retailer's action $a \in \mathbf{A}$.  According to our
assumptions, let us agree that having chosen the action
$a \in \mathbf{A}$ at the current battery level $P_{t}$, the next
level $P_{t+1}$ is modeled as
$$
P_{t+1}= \arg\min_{p\in\mathbf{P}} |p - (P_t + l(a) - \varepsilon_t)|,
$$
from which the transition probabilities in storage levels are induced
by the action $a$ are
\begin{equation} \label{transprob}
\alpha^a_{p,p'} = 
\begin{cases}
	\mathcal N(p+l(a), \varsigma)([p'-\Delta/2, p'+\Delta/2]) & \mbox{if $\underline{p}< p'< \overline{p} $},\\
	\mathcal N(p + l(a), \varsigma)(]-\infty, \underline{p} + \Delta/2] ) & \mbox{if $p'=\underline{p}$}, \\
	\mathcal N(p + l(a), \varsigma)( [ \overline{p} -\Delta/2, +\infty[ ) & \mbox{if $p'=\overline{p}$}
\end{cases}
\end{equation}
where $\mathcal N(p+l(a), \varsigma)$ stands for the probability
measure associated with the normal distribution with mean $p+l(a)$ and
variance $\varsigma^2$.  In a similar manner, the expected excess
$\underline e^a_{p}$ and shortage $\overline e^a_{p}$ of the imbalance
energy can be written as
\begin{equation} \label{excess-shortage}
\begin{array}{c}
\underline e^a_{p} = \int^{\infty}_{\overline{p} +\Delta/2} (x - \overline{p}) \mathcal N(p + l(a), \varsigma)(dx), \\
\overline e^a_{p} = \int_{-\infty}^{\underline{p} - \Delta/2} (\underline{p} - x) \mathcal N(p + l(a), \varsigma)(dx).
\end{array}
\end{equation}
With these definitions, the reward functions are given as above in
(\ref{rewfct}). More specifically, having introduced the rewards
$$
(p, \Pi_t, a) \mapsto  - l(a) \Pi_t - \overline e^a_{p} \overline \Pi_{t}  + 
\underline e^a_{p} \underline \Pi_{t}, \qquad p \in \mathbf{P}, a \in \mathbf{A}, t=0, \dots, T-1
$$
and assuming that at maturity date $T$ the entire energy from storage
capacity can be sold at the forward market, the scrap value is given
by
$$ (p, \Pi_T) \mapsto p \Pi_T, \quad p \in \mathbf{P}.$$
Finally, let us assume that the buy/sell grid prices are constant and
deterministic
$$ 
\underline \Pi_{t}=\underline \Pi,
\quad
\overline \Pi_{t}=\overline \Pi , 
\quad t=0, \dots, T, \quad \text{with} \enspace  0<\underline \Pi < \overline \Pi
$$
to define the reward functions by
\begin{equation}
r_t(p, (z^{(1)}, z^{(2)}), a)=- l(a)  f(t, z^{(2)}) - \overline e^a_{p} \overline \Pi+ \underline e^a_{p} \underline \Pi, \quad
t=0, \dots, T-1, \label{rews}
\end{equation}
for all $a \in \mathbf{A}$, $p \in \mathbf{P}$ and
$(z^{(1)}, z^{(2)}) \in \rr^2$. Finally, introduce the scrap value by
\begin{equation}
r_T(p, (z^{(1)}, z^{(2)}))= p f(T,  z^{(2)}), \quad p \in \mathbf{P}, \quad (z^{(1)}, z^{(2)}) \in \rr^2. \label{sc}
\end{equation}
Note that with the definitions (\ref{rews}), (\ref{sc}) and
(\ref{rec}) our problem of battery storage control is uniquely
determined, whose numerical solution is demonstrated in the next
subsection.

\section{A Case Study} \label{sec:casestudy}
\begin{figure}[h!]
	\includegraphics[height=2.7in, width=0.5\textwidth]{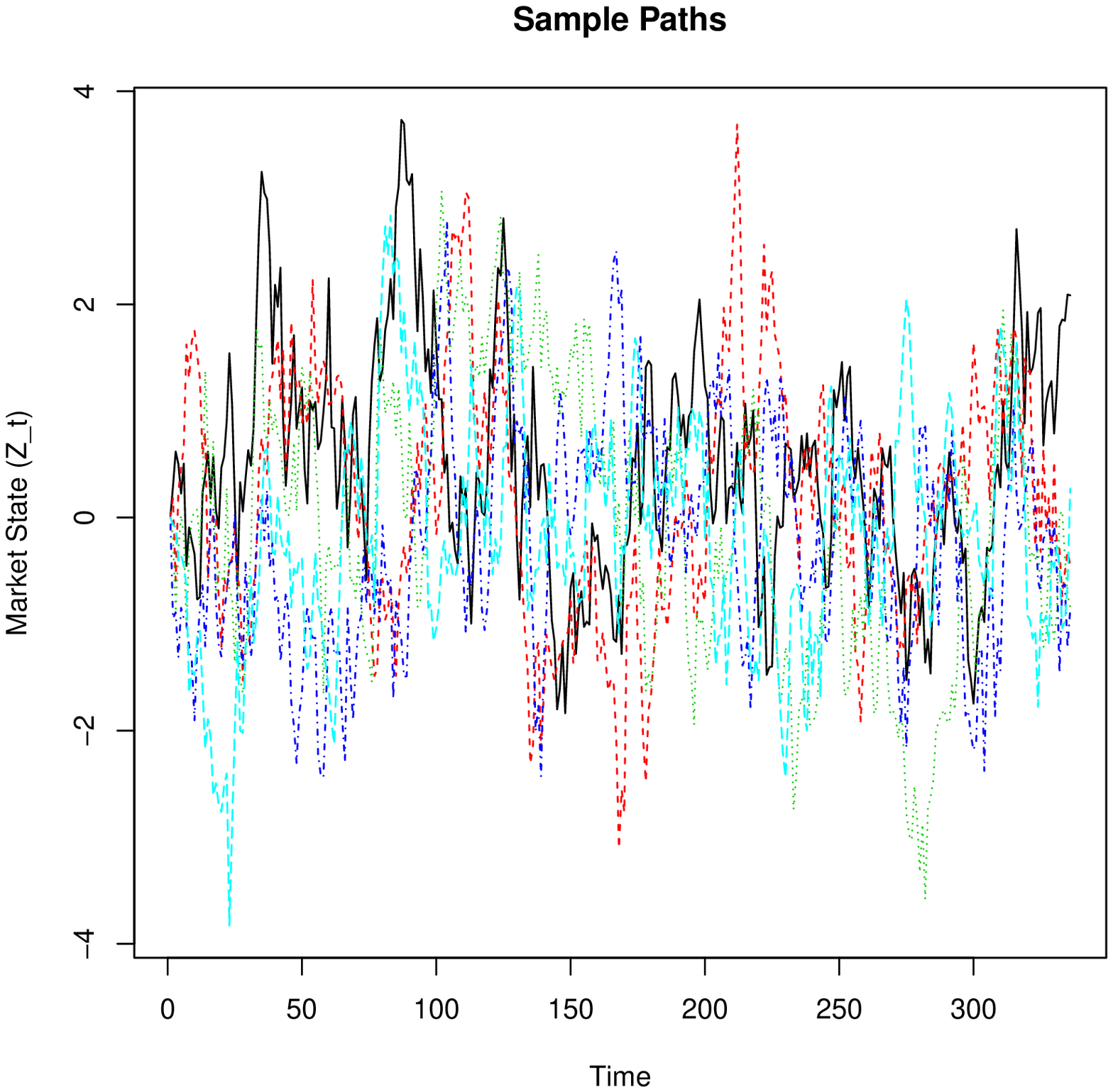}
	\includegraphics[height=2.7in, width=0.5\textwidth]{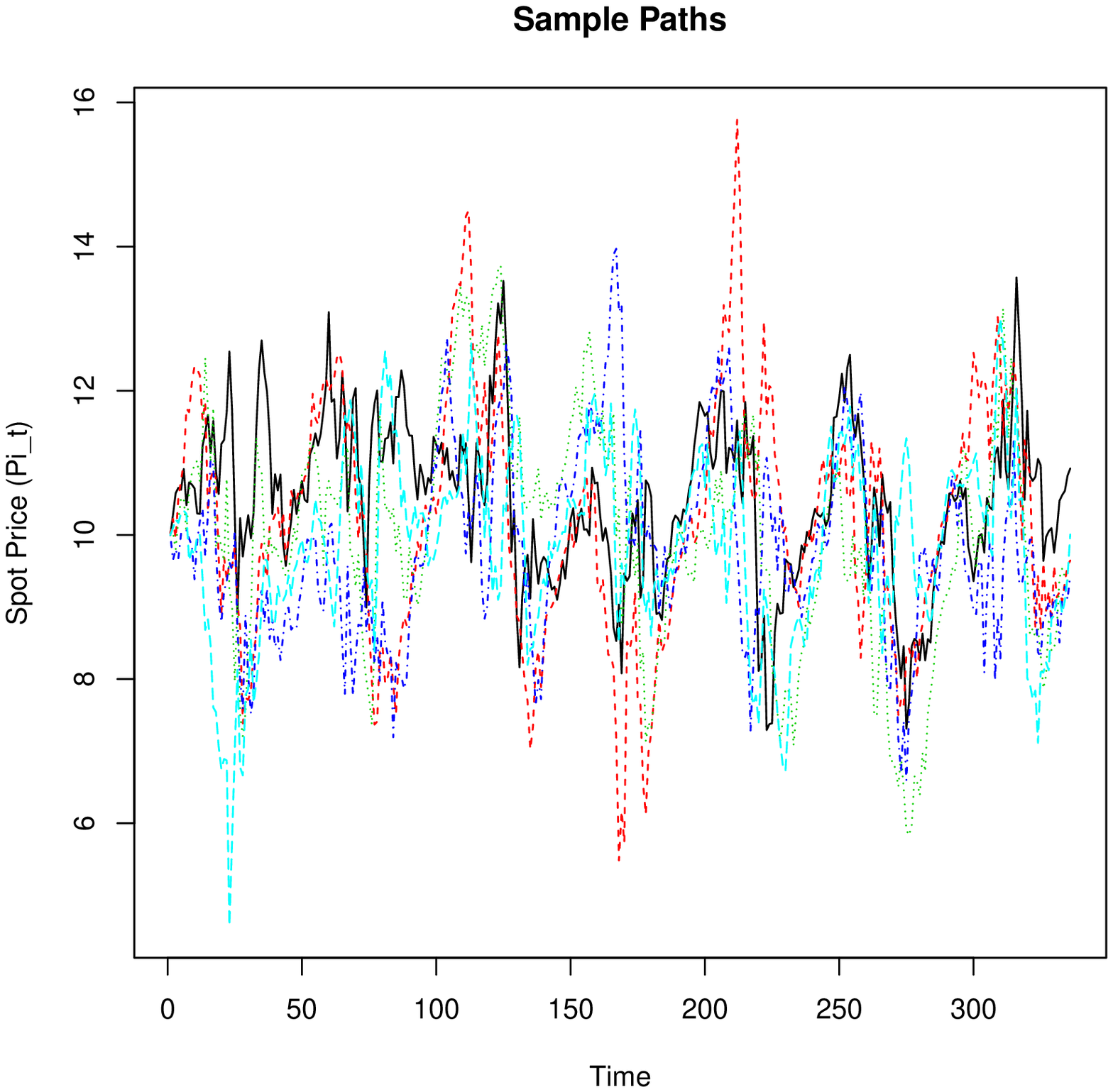}
	\caption{Paths 
          the states $(Z^{(2)}_t)_{t=0}^T$ (left) and the  prices
          $(\Pi_t=f(t, Z^{(2)}_t))_{t=0}^T$, (right). \label{sample_paths}}
\end{figure}

Let us suppose that the battery level are equidistantly discretized with
$\underline p = 0$ MWh, $\overline p = 100$ MWh, and $\Delta = 5$ MWh.
Furthermore, we assume that energy retailer chooses actions
$a \in \mathbf{A}= \{1, 2, \dots, 11\}$ with corresponding safety
margins $l(a) = 5(a - 1)$ MWh for all $a \in \mathbf{A}$.  Further,
consider the time horizon of a week at half-hourly frequency i.e
$t=0, 1, \dots,335=T$ and define the auto-regressive state dynamics as
above, with by $\mu = 0$, $\sigma = 0.5$ and $\phi = 0.9$.  To describe
seasonality, we assume that the affine linear functions are given by
\begin{equation}
f(t, z^{(2)})=  u_t+ v_t z^{(2)}, \qquad t=0, \dots, T, \quad z^{(2)} \in \rr
\label{func}
\end{equation}
with deterministic coefficients
$u_t= 10 + \cos(\frac{2 \pi}{48}t + \frac{3 \pi}{2})$,
$v_t= 1 + \sin(\frac{2 \pi}{48}t + \frac{3 \pi}{2})/2$ for
$t=0, \dots, T$.  Figure \ref{sample_paths} depicts trajectories of
the corresponding price evolution having started $Z^{(2)}_0 = z^{(2)}_0 = 0$.

\begin{figure}[h!]
	\includegraphics[height=2.7in, width=0.5\textwidth]{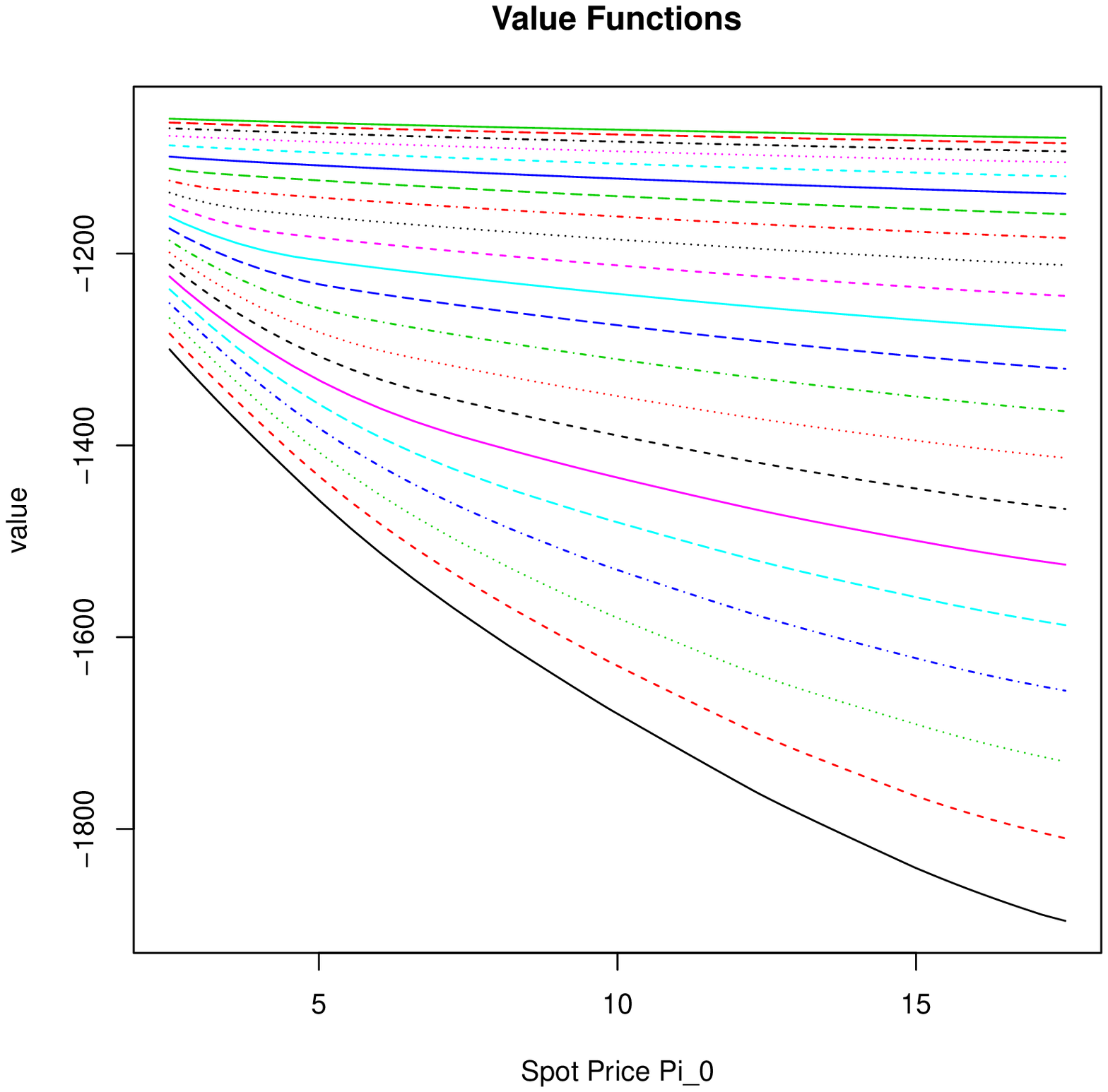}
	\includegraphics[height=2.7in, width=0.5\textwidth]{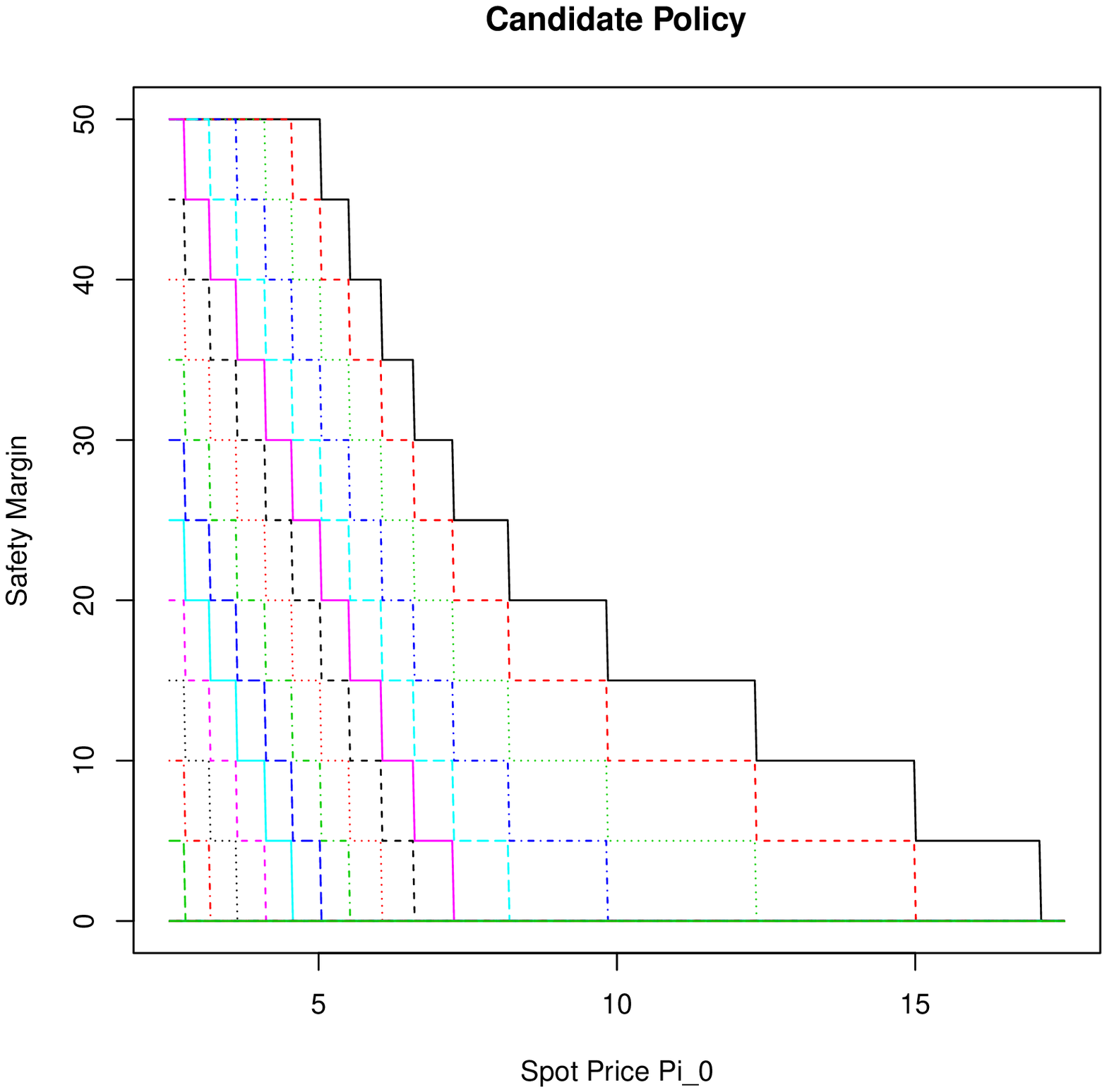}
	\caption{Value functions and optimal policy at $t=0$. Different colours
		for different positions. At each price, the value increases with 
		the current storage level. \label{plots}}
\end{figure}

\begin{table}[h!]
  \centering
  \setlength\tabcolsep{13pt}
  \renewcommand\arraystretch{1}
  \caption{Solution diagnostics for $z^{(2)}=0$ \label{table}. Standard errors in paranthesis.}
  \begin{tabular}{|lrr|}
    \hline
    Battery Level (MWh) & Lower Bound & Upper Bound \\
    \hline
   0  &-1679.759 (0.042) & -1679.756 (0.042) \\
   5  &-1629.759 (0.042) & -1629.756 (0.042) \\
  10  &-1579.759 (0.042) & -1579.756 (0.042) \\
  15  &-1529.759 (0.042) & -1529.756 (0.042) \\
  20  &-1480.069 (0.042) & -1480.066 (0.042) \\
  25  &-1433.475 (0.041) & -1433.472 (0.041) \\
  30  &-1389.587 (0.041) & -1389.583 (0.041) \\
  35  &-1348.411 (0.041) & -1348.408 (0.041) \\
  40  &-1310.032 (0.041) & -1310.028 (0.041) \\
  45  &-1274.505 (0.041) & -1274.502 (0.041) \\
  50  &-1241.857 (0.040) & -1241.853 (0.040) \\
  55  &-1212.091 (0.040) & -1212.088 (0.040) \\
  60  &-1185.201 (0.040) & -1185.197 (0.040) \\
  65  &-1161.168 (0.040) & -1161.165 (0.040) \\
  70  &-1139.971 (0.039) & -1139.968 (0.039) \\
  75  &-1121.586 (0.039) & -1121.583 (0.039) \\
  80  &-1105.989 (0.039) & -1105.986 (0.039) \\
  85  &-1093.160 (0.039) & -1093.157 (0.039) \\
  90  &-1083.071 (0.039) & -1083.068 (0.039) \\
  95  &-1075.638 (0.039) & -1075.634 (0.039) \\
 100  &-1070.639 (0.039) & -1070.636 (0.039) \\
    \hline
  \end{tabular}
\end{table}
Assume that the standard deviation of the consumer's demand prediction
error is $\varsigma = 10$ and define grid prices by
$\overline \Pi = 20$ and $\underline \Pi = 0$, respectively. 
With these quantities, we apply Algorithm \ref{algo1} to a state space
grid $G$ containing $501$ points equally distributed on the line
connecting the points $(1, -15)$ and $(1, 15)$.  Furthermore, we
discretize the Normal distribution in the disturbance matrices in
terms of $10000$ equidistant quantiles.  The solution diagnostics in
Table \ref{table} is generated by Algorithm \ref{algo2} with $100$
sample paths for the price and $100$ subsimulations for each path at
each non-terminal decision epoch.  All results are obtained using
authors' R package \cite{rcss}.  We obtained tight bounds and low
standard errors, which certify that our results are sufficiently close
to the true solution. Figure \ref{plots} illustrates the policy values
and the decision structure of our nearly-optimal policy.

\begin{figure}[h!]
	\includegraphics[height=2.7in,width=0.5\textwidth]{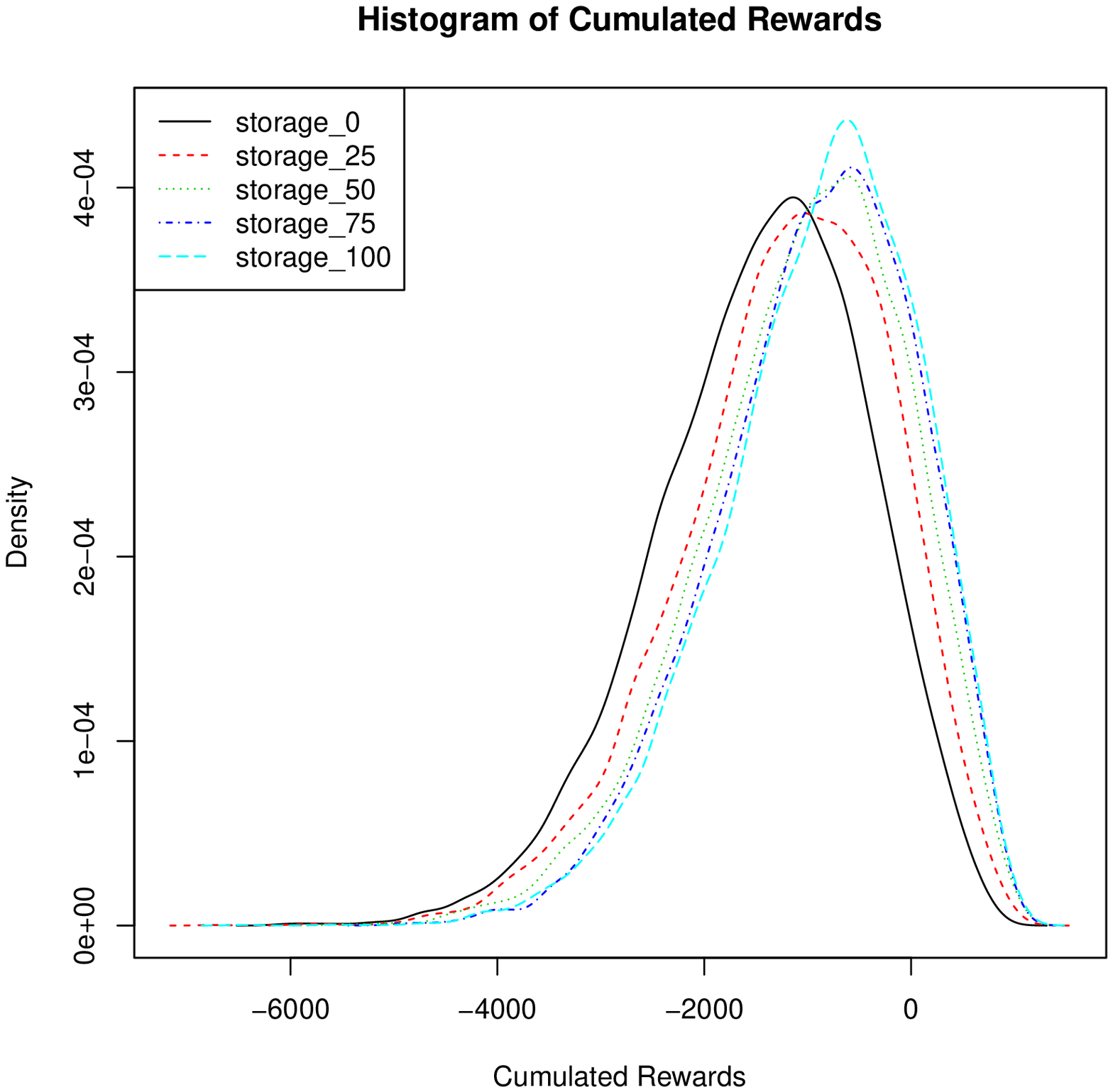}
	\includegraphics[height=2.7in,width=0.5\textwidth]{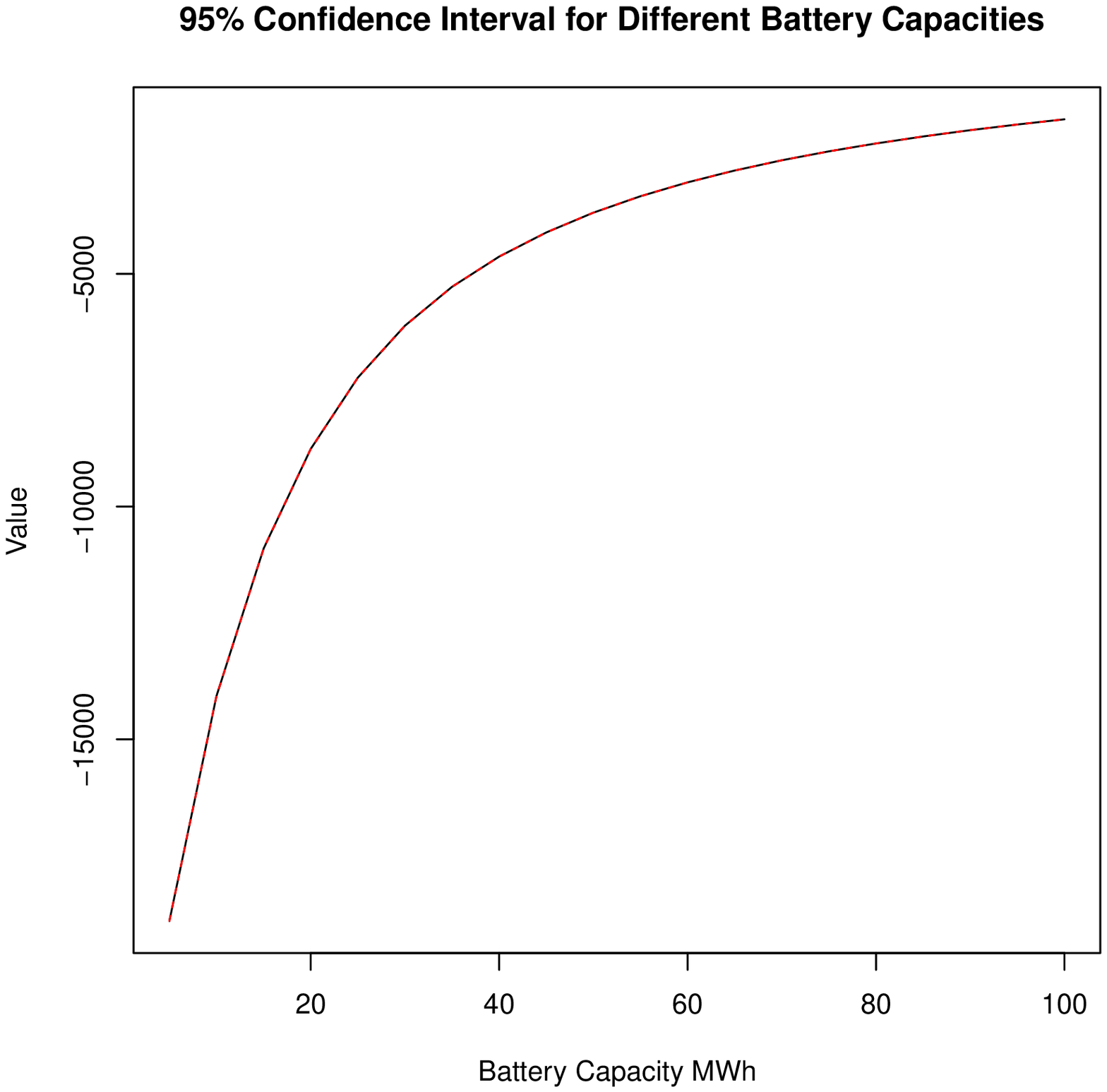}
	\caption{The left plot shows the histogram of the cumulated
          rewards on 10000 sample paths for different starting levels
          using the prescribed policy. The right depicts two curves,
          enclosing the 95\% confidence region for the cumulated
          reward, depending on storage capacity.\label{plot_hist}}
\end{figure}

Finally, let us discuss a typical economic application. Having solved
the problem of optimal operational management, the questions of
investment and capacity allocation can be addressed.  Here, the first
step is an estimation of the random revenue profile, followed in the
second step by a detailed risk analysis from potential side effects of
operational flexibility on the entire portfolio of existing physical
and financial assets. Thereby, the hedging value of flexibility (see
\cite{Doege2009936}) is essential.  In our study, we illustrate this
first step using Monte-Carlo simulations.  Having calculated the
revenue of our approximately optimal strategy on 10000 randomly
generated path scenarios, its density histogram is plotted in the left
plot of Figure \ref{plot_hist} for different starting storage levels.
In line with our expectations, we observe in this figure that a higher
initial battery level yields a higher cumulated reward. This is also
seen from the left plot of the Figure \ref{plots} which shows that the
value function at any price is increasing in storage level (which
correspond to different curves). The right graph of Figure \ref{plots}
indicates that a higher initial storage level also yields a lower
safety margin in the optimal strategy at the initial time.
\\
\\
{\bf Remark:}\footnote{The authors thank an anonymous referee for
  pointing to this analysis.} The right plot in Figure \ref{plot_hist}
shows the dependence of the cumulated rewards on the battery size. In
this graph, we gradually increase the capacity from 5 MWh to 100MWh
and determine the value function for an initially empty storage (for
$\phi=0.9$, $z^{(2)}_0=0$).  The concave curve depicted in this graph
shows that the value grows with the capacity at a rate which is
steadily decreasing.  Such insight may be very valuable in
practice. For instance, the optimal size for battery deployment would
typically result from equating the marginal value of the storage to
its marginal cost.  Notice however, that the true practical value of
such optimization usually stems from risk hedging effects within
agent's energy generation portfolio, whose analysis now becomes
possible, using a reliable strategy optimization provided by our
concepts.

\subsection{Optimal Trading and Storage} \label{sec:discussion}

In this section, we discuss the impact of storage facility on optimal
energy trading. Recall that in our model, the agent indirectly
controls the storage level in terms of energy trading. The point is
that the storage absorbs some (if not all) of any unexpected demand
which adds some flexibility in the price choice when energy is
purchased.  To this end, we investigate the impact of battery size,
comparing a small (5MWh) against large (100MWh) storage.

\begin{table}[h!]
  \centering
  \setlength\tabcolsep{3pt}
  \caption{Bounds on  cumulated rewards  estimation for  different parameters.
   \label{table_ar}}
  \begin{tabular}{l|cc|cc}
    & \multicolumn{2}{c}{Small Battery 5 MWh}&\multicolumn{2}{c}{Large Battery 100 MWh}\\
    $\phi$ & Lower & Upper & Lower & Upper \\
    \hline
    0.9 &  -18904.06 (0.151) & -18904.06 (0.151) & -1679.759 (0.042) & -1679.756 (0.042) \\
    0.6 &  -19004.19 (0.073) & -19004.19 (0.073) & -1682.616 (0.037) & -1682.609 (0.037) \\
    0.3 &  -19017.53 (0.060) & -19017.52 (0.059) & -1679.807 (0.038) & -1679.799 (0.039) \\
    0.1 &  -19019.21 (0.057) & -19019.20 (0.057) & -1676.744 (0.042) & -1676.732 (0.042) \\
  \end{tabular}
\end{table}

\begin{figure}[h!]
	\includegraphics[height=3.2in,width=\textwidth]{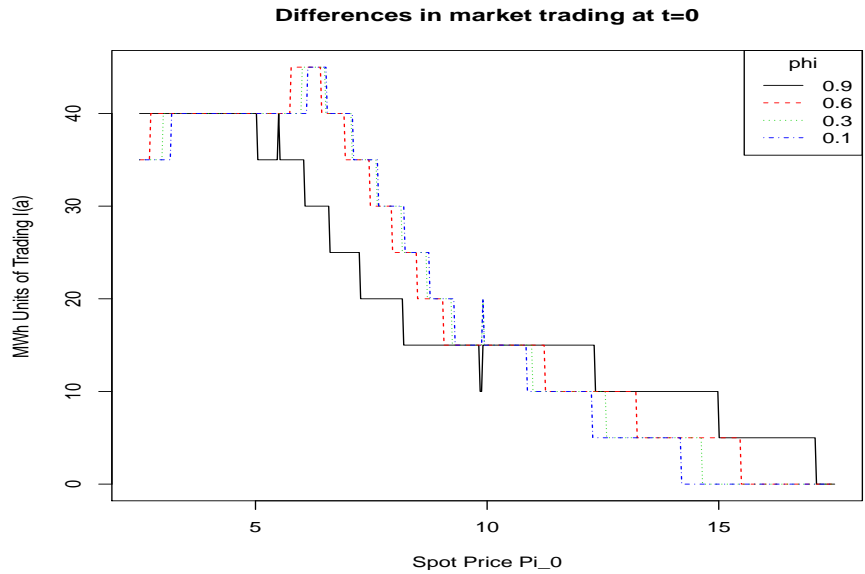}
	\caption{Differences in optimal trading at $t=0$. Both
          batteries intially empty. \label{plot_startpolicy}}
\end{figure}

In our state dynamics (\ref{rec}), the parameter $\phi$ controls the
speed of mean reversion for the energy price. Thereby, the lower
levels of $\phi$ lead to stronger mean reversion with more frequent
return to the seasonal price component of the price. Table
\ref{table_ar} compares the expected cumulated rewards for a small
battery with 5 MWh capacity against a large battery of size 100 MWh
under the assumptions that both batteries are initially empty. This
table shows that there is a very substantial benefit of extra storage
capacity for all levels of mean reversion. However, mean reversion
seems to have very little impact on the expected cumulated rewards.
Further, Figure \ref{plot_startpolicy} illustrates the role of storage
capacity for energy purchase.  Here, we depict the difference between
the safety margins optimally entered at time $t=0$ for empty storage.
We observe that having a large battery capacity allows to buy more
energy in advance. Note also that this effect decreases with
increasing price. Thereby, the impact of mean reversion parameter is
low again.

\begin{figure}[h!]
	\includegraphics[height=2.7in, width=0.5\textwidth]{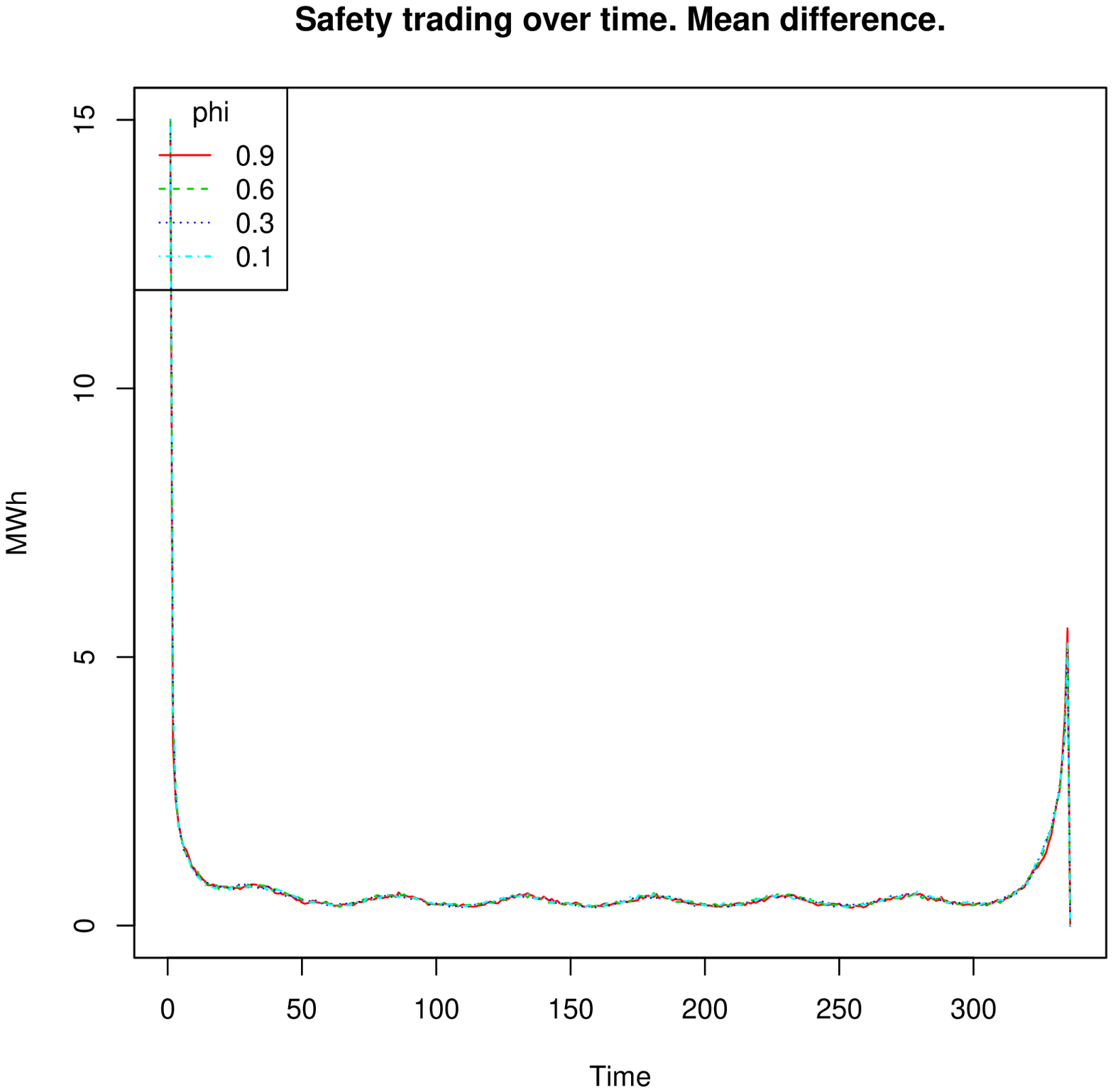}
	\includegraphics[height=2.7in, width=0.5\textwidth]{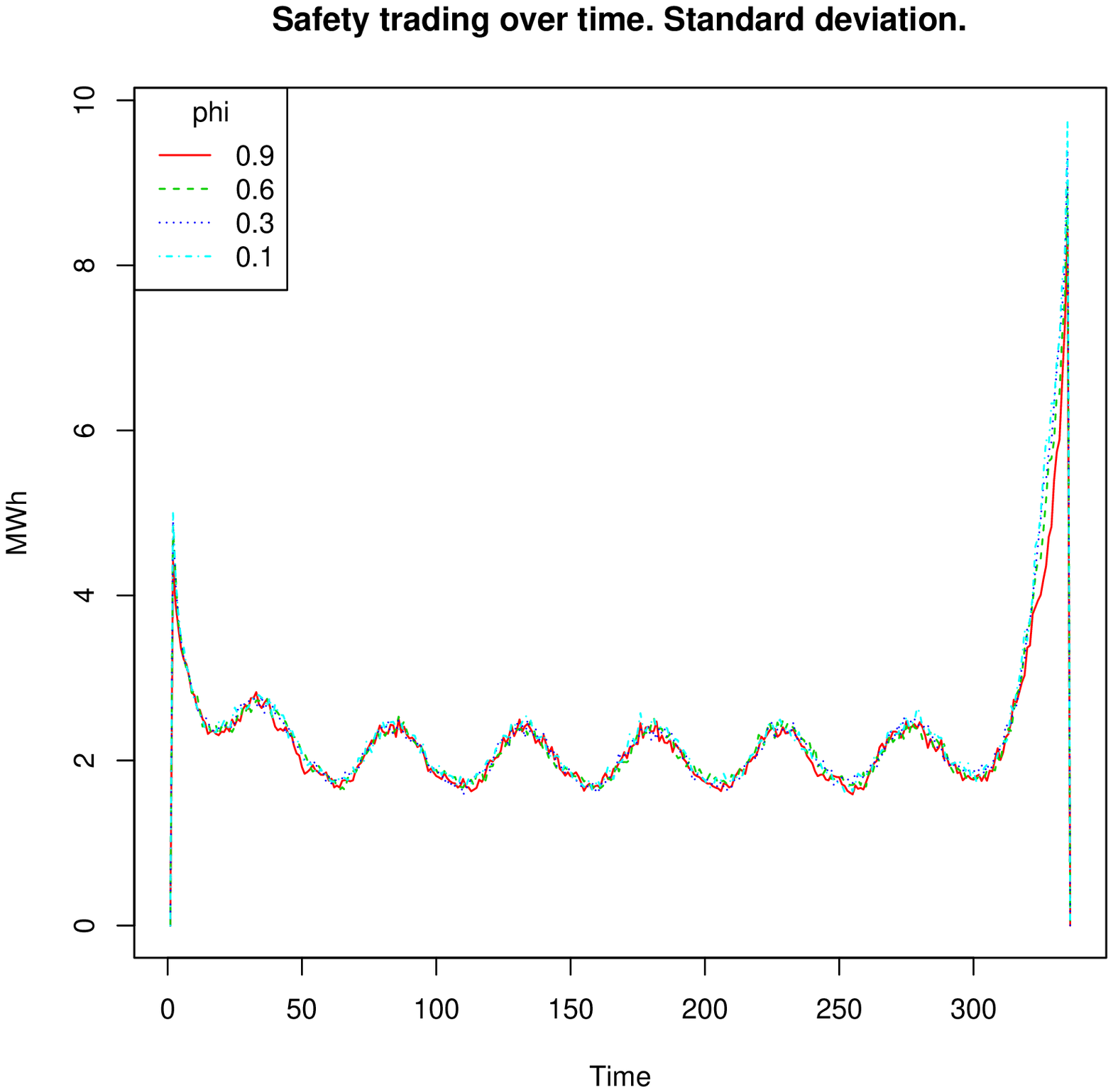}
	\caption{Averaged difference (10000 scenarios)  in safety margin between the large
          and small battery (both initially empty) over time for
          different mean reversions.  The empirical mean is  on
          the left and the empirical  standard deviation is on the
          right. \label{plot_action}}
\end{figure}

\begin{figure}[h!]
	\includegraphics[height=2.7in, width=0.5\textwidth]{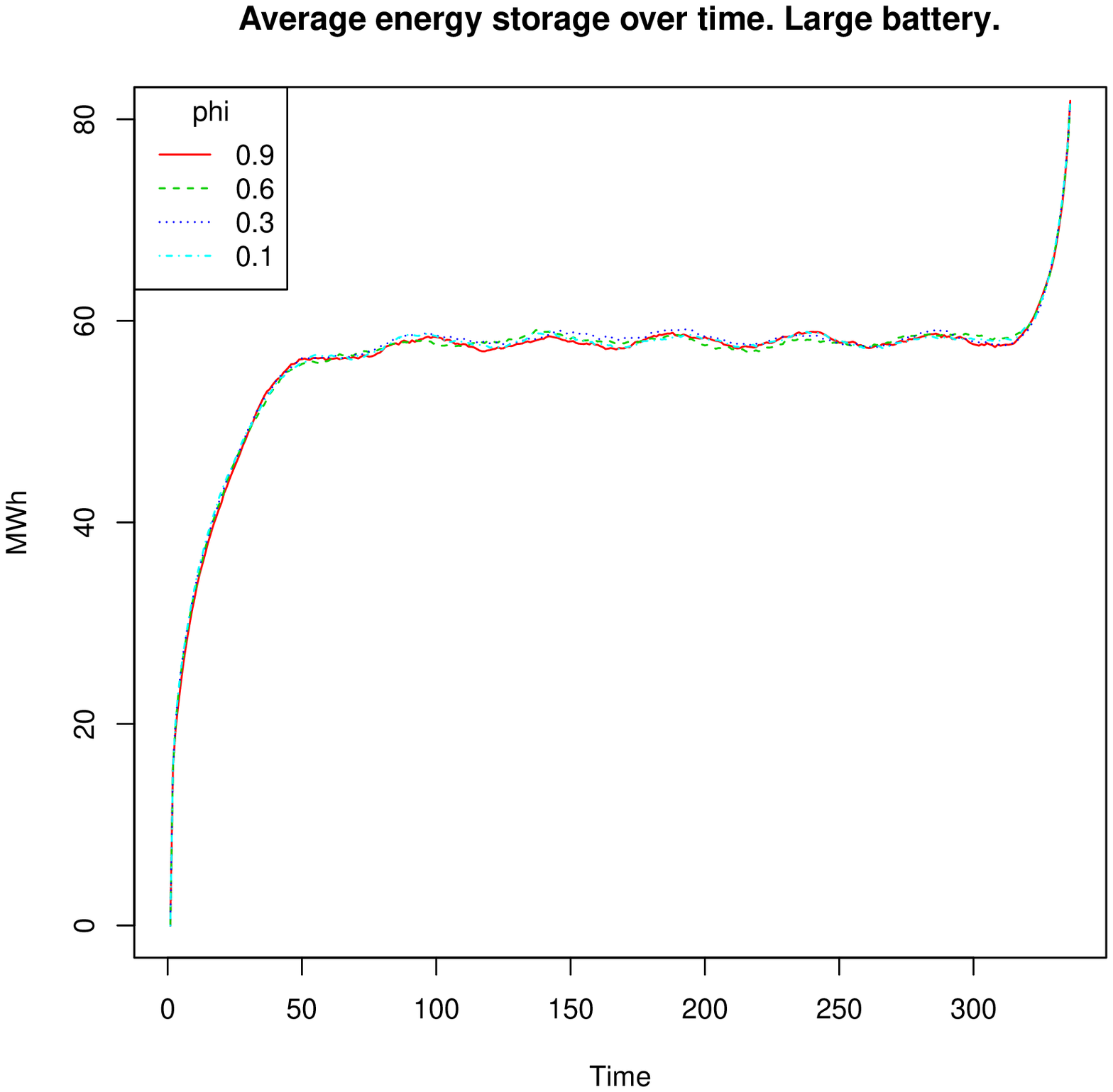}
	\includegraphics[height=2.7in, width=0.5\textwidth]{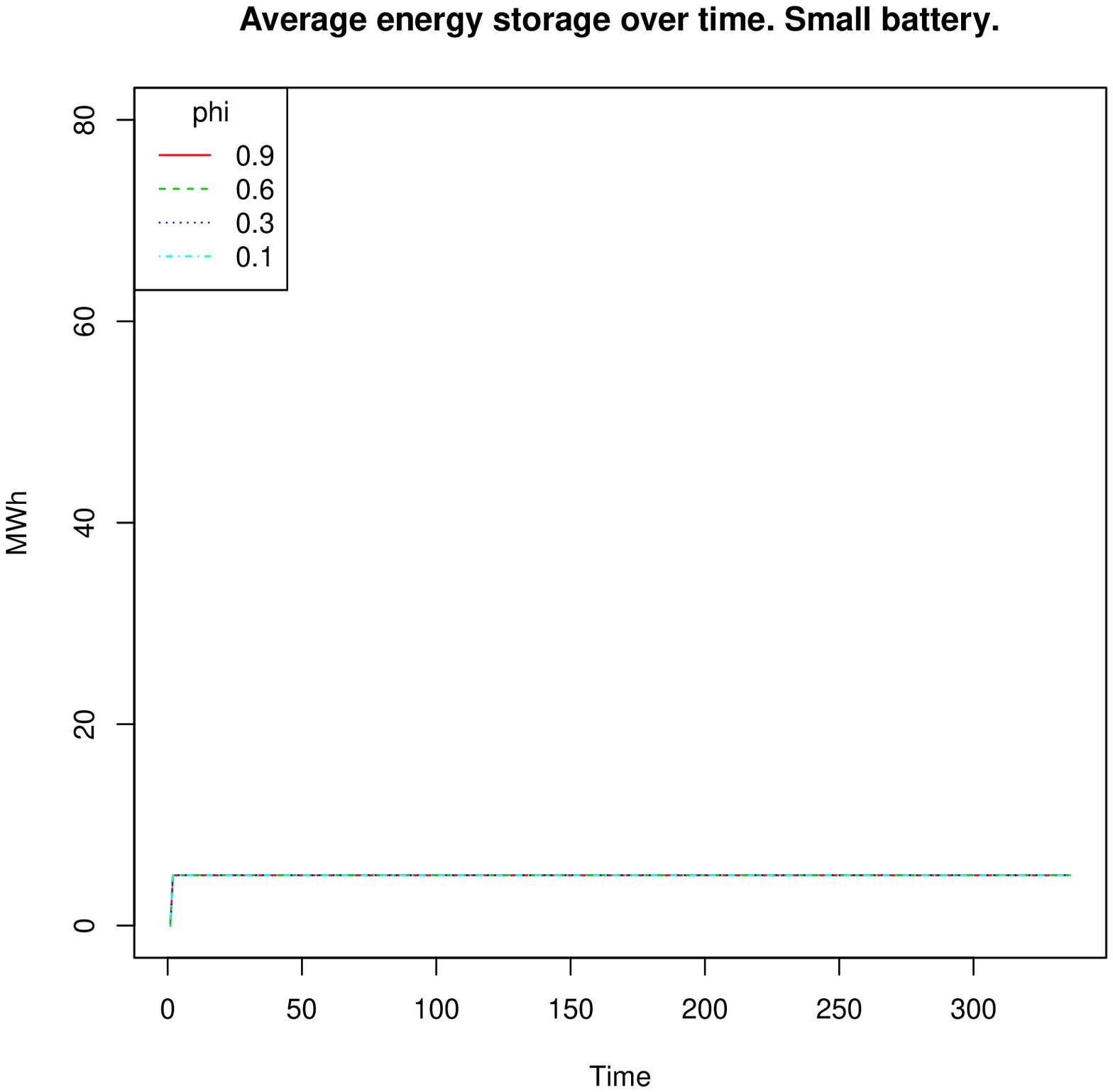}
	\caption{The averaged evolution (10000 scenarios) of the energy storage level for the 
          large battery, on the left, and for  small battery on the right.
           \label{plot_position}}
\end{figure}

Finally, we examine in Figure \ref{plot_action} the averaged behavior
(on $10000$ scenarios) of safety margins in dependence on storage
capacity for different mean reversion parameters.  In line with the
previous observations, the speed of mean reversion seems to have only
little impact.  Remarkably, there is a clear seasonal pattern in the
difference of optimal safety margins.  This is caused by the seasonal
nature of prices as shown in Figure \ref{sample_paths}. Namely, for
large storage capacity, the optimal trading follows price seasonality
stronger than if the storage is small.  This issue is also obvious
from Figure \ref{plot_position} which shows that for large capacity it
is optimal to keep a certain intermediate level whereas small storage
must be filled right at the beginning to hedge against unexpected
demand fluctuations.  We also observe that close to maturity there is
an attempt to fill the storage in order to benefit from price at the
end of the time period. 

\begin{table}[h!]
  \centering
  \setlength\tabcolsep{1pt}
  \caption{Bounds on cumulated rewards under $\phi = 0.9$ and $z^{(2)} = 0$.
   \label{table_arbitrage}}
  \begin{tabular}{l|cc|cc}
    Capa- & \multicolumn{2}{c}{Non-zero scrap value}&\multicolumn{2}{c}{Zero scrap value}\\
    city & Lower & Upper & Lower & Upper \\
    \hline
 10 & -14068.958 (0.115) & -14068.957 (0.115) & -14124.612 (0.115) & -14124.611 (0.115) \\
 20 &  -8762.276 (0.078) &  -8762.275 (0.077) &  -8879.116 (0.078) &  -8879.115 (0.078) \\
 30 &  -6114.388 (0.049) &  -6114.388 (0.049) &  -6292.050 (0.049) &  -6292.049 (0.049) \\
 40 &  -4629.497 (0.039) &  -4629.496 (0.039) &  -4866.371 (0.039) &  -4866.370 (0.039) \\
 50 &  -3685.724 (0.033) &  -3685.723 (0.033) &  -3980.018 (0.033) &  -3980.017 (0.033) \\
 60 &  -3033.977 (0.030) &  -3033.977 (0.030) &  -3384.379 (0.029) &  -3384.379 (0.029) \\
 70 &  -2559.781 (0.028) &  -2559.781 (0.028) &  -2965.728 (0.027) &  -2965.728 (0.027) \\
 80 &  -2198.558 (0.031) &  -2198.557 (0.031) &  -2660.035 (0.027) &  -2660.034 (0.027) \\
 90 &  -1912.817 (0.035) &  -1912.815 (0.035) &  -2430.169 (0.029) &  -2430.168 (0.029) \\
100 &  -1679.759 (0.042) &  -1679.756 (0.042) &  -2253.495 (0.033) &  -2253.493 (0.033) \\
  \end{tabular}
\end{table}

\noindent {\bf Remark:}\footnote{The authors thank an anonymous referee for raising this issue.} 
In our finite horizon setting, three phases are observable from Figure
\ref{plot_position}. After an initial charge to an ''optimal
intermediate level'', the battery is used to absorb unexpected
demand/supply fluctuations while remaining close to this
level. However, closer to maturity, the storage is further filled to
take advantage of the price at the end of the time horizon.  On this
account, it is important to compare the portion of the cumulative
reward which results from intermediate balancing to that earned from
selling energy at maturity. To investigate this problem, alter the
scrap value definition (\ref{sc}) to
\begin{equation}
r_T(p, (z^{(1)}, z^{(2)}))= 0 \quad p \in \mathbf{P}, \quad (z^{(1)}, z^{(2)}) \in \rr^2. \label{sc1}
\end{equation}
With this change, any energy remaining at the end is worthless. Table
\ref{table_arbitrage} compares the expected cumulative rewards of the
original (\ref{sc}) problem to that with (\ref{sc1}) for different
storage capacities. It is not surprising that a larger battery allows
exploiting the remaining energy to a greater extent. However comparing
both columns in this table, we observe that a very significant part of
the battery value results from the energy balancing.

\section{Future Research and Conclusion} \label{sec:conclusion} 

Electrical storages have the potential to essentially change the
nature of electricity trading and may have profound impact on energy
price dynamics.  This paper provides quantitative concepts to better
understand and analyze this development.  We demonstrate that using
our algorithmic approach to battery storage management, a detailed and
accurate strategy optimization is possible.  Further details, such as
modelling uncertainties in grid prices, costs of deep discharge
affecting battery's life time, and stochastic futures price dynamics
can be incorporated. The authors will address these exciting topics in
future research.

\
\section{References}


\bibliographystyle{plain}
\bibliography{co2-3-2,battery}

\newpage
\appendix
\section{R Script for Table \ref{table}}  \label{appendix} 
\noindent
The following code was used for Table \ref{table}. On Linux Ubuntu
16.04 with Intel i5-5300U CPU @2.30GHz and 16GB of RAM, the script
below takes less than 20 seconds to run and requires the installation of
the \emph{'rcss'} package \cite{rcss}.
\\
\\
The first part of the script returns the value function approximations
using Algorithm \ref{algo1}. It takes roughly 5 seconds to run.
\begin{lstlisting}
## Remove existing objects and load R package
rm(list = ls()); gc(); library(rcss)
## Grid
grid <- cbind(rep(1, 501), seq(-15, 15, length = 501))
## Battery
battery <- seq(0, 100, by = 5)
## Standard deviation for the consumer demand
std <- 10
## Safety margins
safety <- seq(0, 50, length = 11)  ## safety
## Transition probabilities for controlled Markov chain
control <- array(data = 0, dim = c(21, 11, 21))
for (p in 1:21) {
    for (a in 1:11) {
        temp <- battery[p] + safety[a]  ## center of normal distribution
        control[p,a,1] <- pnorm(0 + 5/2, temp, std)
        control[p,a,21] <- 1 - pnorm(100 - 5/2, temp, std)
        for (pp in 2:(21-1)) {
            control[p,a,pp] <- pnorm(battery[pp] + 5/2, temp, std) -
                pnorm(battery[pp] - 5/2, temp, std)
        }
    }
}   
## Functions to calculate expected excess and shortage energy demand
erf <- function(x){  ## error function
    return(2 * pnorm(x * sqrt(2)) - 1)
}
Excess <- function(pos, act) {
    temp1 <- 100 + 5/2
    temp2 <- pos + act
    result <- std/sqrt(2*pi) * exp(-(temp1-temp2)^2/(2*std^2)) +
        (temp2 - 100)/2 * (1 - erf(1/sqrt(2*std^2) * (temp1 - temp2)))
    return(result)
}
Shortage <- function(pos, act) {
    temp1 <- 0 - 5/2
    temp2 <- pos + act
    result <- std/sqrt(2*pi) * exp(-(temp1-temp2)^2/(2*std^2)) +
        (0 - temp2)/2 * (erf(1/sqrt(2*std^2) * (temp1 - temp2)) + 1)
    return(result)
}
## Expected excess and shortage energy demand
excess <- matrix(data = NA, nrow = 21, ncol = 11)
shortage <- matrix(data = NA, nrow = 21, ncol = 11)
for (p in 1:21) {
    for (a in 1:11) {
        excess[p,a] <- Excess(battery[p], safety[a])
        shortage[p,a] <- Shortage(battery[p], safety[a])
    }
}
## Subgradient representation of reward functions
u_t <- 10 + cos((0:335) * 2*pi/48 + 3*pi/2)
v_t <- 1 + (sin((0:335) * 2*pi/48 + 3*pi/2))/2
reward <- array(0, dim = c(501, 2, 21, 11, 336))
for (p in 1:21) {
    for (a in 1:11) {
        for (t in 1:335) {
            reward[,1,p,a,t] <- -safety[a] * u_t[t] - shortage[p, a] * 20
            reward[,2,p,a,t] <- -safety[a] * v_t[t]
        }
    }
    ## Scrap reward
    reward[,1,p,,336] <- battery[p] * u_t[336]
    reward[,2,p,,336] <- battery[p] * v_t[336]
}
## Parameters for AR(1) process (Z_t)
mu <- 0
sigma <- 0.5
phi <- 0.9
## Disturbances (W_t)
disturb_weight <- rep(1 / 10000, 10000)  ## probability weights
disturb <- array(matrix(c(1, 0, 0, phi), ncol = 2, byrow = TRUE), dim = c(2, 2, 10000))
quantile <- qnorm(seq(0, 1, length = (10000 + 2))[c(-1, -(10000 + 2))])
disturb[2, 1,] <- mu + sigma * quantile
r_index <- matrix(c(2, 1), ncol = 2)  ## randomness index
## Fast bellman recursion
bellman <- FastBellman(grid, reward, control, disturb, disturb_weight, r_index)
\end{lstlisting}

\vspace{3mm}
	
\noindent
The second part of the script computes the lower and upper bound
estimates according to Algorithm \ref{algo2}. It takes roughly 10 seconds 
to run.
\begin{lstlisting}
## Exact reward function
Reward <- function(state, time) {
    output <- array(0, dim = c(nrow(state), 21 * 11))
    if (time == 336) {
        for (p in 1:21) {
            output[,(p-1) * 11 + (1:11)] <-
                battery[p] * (u_t[time] + v_t[time] * state[,2])
        }
        return(output)
    }
    for (p in 1:21) {
        for (a in 1:11) {
            output[,(p-1) * 11 + a] <-
                -safety[a] * (u_t[time] + v_t[time] * state[,2]) -
                shortage[p,a] * 20
        }
    }
    return(output)
}
## Generate sample path disturbances
set.seed(12345)
path_disturb <- array(matrix(c(1, 0, 0, phi), ncol = 2, byrow = TRUE),
                      dim = c(2, 2, 335, 100))
rand <- rnorm(335 * 100 / 2)
rand <- c(rand, -rand)
path_disturb[2, 1,,] <- mu + sigma * rand
## Specifying subsimulation disturbances
subsim_weight <- rep(1 / 100, 100)
subsim_disturb <- array(matrix(c(1, 0, 0, phi), ncol = 2, byrow = TRUE),
                        dim = c(2, 2, 100, 100, 335))
rand <- rnorm(100 * 100 * 335 / 2)
rand <- as.vector(rbind(rand, -rand))
subsim_disturb[2, 1,,,] <- mu + sigma * rand
## Generate sample paths for uncontrolled process
start <- c(1, 0)  ## z_0
path <- Path(start, path_disturb)
path_nn <- Neighbour(matrix(path, ncol = 2), grid, 1, "kdtree", 0, 1)$indices
## Candidate policy for sample paths
path_action <- PathPolicy(path, path_nn, control, Reward, bellman$expected, grid)
## Computing martingale increments
time2 <- proc.time()
mart <- FastMartingale(bellman$value, path, path_nn, subsim_disturb,
                       subsim_weight, grid, control = control)
## Calculating the primal and dual values
duality <- Duality(path, control, Reward, mart, path_action)
time2 <- proc.time() - time2
## Printing the results
diagnostics <- matrix(data = NA, nrow = 21, ncol = 4)
for (p in 1:21) {
    diagnostics[p, 1] <- mean(duality$primal[1, p,])
    diagnostics[p, 2] <- sd(duality$primal[1, p,])/sqrt(100)
    diagnostics[p, 3] <- mean(duality$dual[1, p,])
    diagnostics[p, 4] <- sd(duality$dual[1, p,])/sqrt(100)
}
print(round(diagnostics,3))
\end{lstlisting}

\end{document}